\documentclass[graybox]{svmult}

\usepackage{mathptmx}       
\usepackage{helvet}         
\usepackage{courier}        
\usepackage{makeidx}         
\usepackage{graphicx}        
\usepackage{multicol}        
\usepackage[bottom]{footmisc}

\usepackage{cite}
\usepackage{fancyvrb}        
\usepackage{siunitx}
\usepackage{amsmath,amssymb,bm} 
\usepackage{subcaption}
\newcommand \bfB {\mathbf{B}}
\newcommand \bfA {\mathbf{A}}
\newcommand \bfx {\mathbf{x}}
\newcommand \dd {{\mathrm d}}
\newcommand \bfc {\mathbf{c}}
\newcommand \Np {N_{\mathrm p}}
\newcommand \Ns {N_{\mathrm s}}
\newcommand \bfxhat {\mathbf{\widehat{x}}}
\newcommand \bfchat {\mathbf{\widehat{c}}}
\newcommand \bfw {\mathbf{w}}
\newcommand \bfy {\mathbf{y}}
\newcommand \xhat {{\widehat{x}}}
\DeclareMathAlphabet{\mathcalOb}{OMS}{cmsy}{b}{n}
\DeclareBoldMathCommand \cA{{\mathcalOb A}}
\DeclareBoldMathCommand \cB{{\mathcalOb B}}
\DeclareBoldMathCommand \cJ{{\mathcalOb J}}
\DeclareBoldMathCommand \cQ{{\mathcalOb Q}}
\DeclareBoldMathCommand \cC{{\mathcalOb C}}
\newcommand \Ts {T_{\mathrm s}}
\newcommand \fs {f_\mathrm{s}}
\newcommand \bfX {\mathbf{X}}

\begin{document}

\title*{Parallel-in-Time Simulation of Power Converters Using Multirate PDEs}
\titlerunning{Parallel-in-Time Simulation of Power Converters Using Multirate PDEs}
\author{Andreas Pels 
\and
Iryna Kulchytska-Ruchka
\and 
Sebastian Sch{\"o}ps 
}
\institute{Andreas Pels, Iryna Kulchytska-Ruchka, Sebastian Sch{\"o}ps \at Computational Electromagnetics Group, Technical University of Darmstadt, Schloßgartenstr. 8, 64289 Darmstadt, Germany,  
\email{\{pels, kulchytska, schoeps\}@temf.tu-darmstadt.de}
 }

%
%
\maketitle

\abstract*{}

\abstract{
This paper presents a numerical algorithm for the simulation of
pulse-width modulated power converters via parallelization in time domain. 
The method applies the multirate partial differential equation approach on the coarse grid of the 
(two-grid) parallel-in-time algorithm Parareal. Performance of the proposed 
approach is illustrated via its application to a DC-DC converter.
}

\section{Introduction}
\label{pels:sec_general}
\begin{figure}[t]
  \begin{subfigure}{0.45\textwidth}
    \centering
    \includegraphics[width=\columnwidth, trim=0 -3em 0 0]{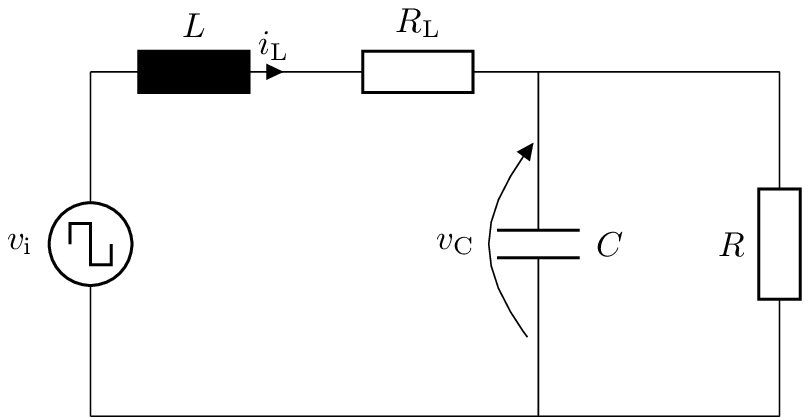}
    \caption{}
    \label{pels:fig:buck_conv}
  \end{subfigure}
  ~~
  \begin{subfigure}{0.55\textwidth}
    \centering
	\includegraphics{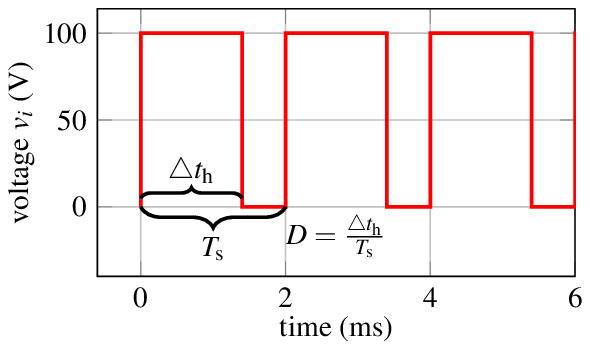}
    \caption{}
    \label{pels:fig:pwmSignal}
  \end{subfigure}
  \caption{Power converter model with pulsed voltage source: (a) Circuit of a simplified buck converter. Transistor switching is modeled as pulsed voltage source. (b) PWM generated pulsed voltage.}
\end{figure}

Switch-mode power converters are devices which convert electric voltages or currents between different levels. For this purpose they use transistors to switch on and off the input voltage or current to obtain the desired average voltage or current at the output of the converter. A technique called pulse-width modulation (PWM) is often utilized to control the transistors, i.e., to generate the pulsed voltage from a given carrier and reference.  An exemplary circuit of a buck converter (DC-DC converter) is depicted in Fig.~\ref{pels:fig:buck_conv} along with its solution in Fig.~\ref{pels:fig:solution}. It consists of fast periodically varying ripples and a slowly varying envelope. The simulation of these power converters with conventional time stepping is computationally expensive since a high number of time steps is necessary to resolve the fast variations induced by the transistor switching.

This paper proposes the simulation of power converters using a combination of two methods, namely the parallel-in-time algorithm Parareal \cite{Lions_2001aa} and a multirate approach based on Multirate Partial Differential Equations (MPDEs) \cite{Pels_2019aa}. This is accomplished via the application of the MPDE approach on the coarse grid of Parareal. It allows the coarse propagator to obtain a more precise solution given the PWM input signal, in contrast to the standard coarse propagator when using a large time step on the original system of equations.

\medskip

The paper is organized as follows: first we introduce our model problem with pulsed excitation in Section~\ref{pels:sec:problem}, then in Section~\ref{pels:sec:parareal} the Parareal method is summarized, Section \ref{pels:sec:mpdes} proposes the usage of MPDEs as coarse propagators for Parareal that can deal with pulsed right-hand sides and finally Section~\ref{pels:sec:numerics} discusses a numerical example before concluding the paper.

\section{Power Converter Model}
\label{pels:sec:problem}
Switch-mode power converters, {which convert AC to DC, DC to AC, AC to AC, or DC to DC voltages, are frequently used devices}. They use power electronic switches to periodically switch the input voltage on and off to regulate the output voltage. For example a buck converter (DC-DC converter) transforms a given voltage to a lower output voltage. It consists of a part that generates a pulsed voltage $v_i$ and a filter circuit. The latter is shown in Fig.~\ref{pels:fig:buck_conv}. The pulsed voltage, see Fig.~\ref{pels:fig:pwmSignal}, is often generated using PWM. Important quantities defining the pulsed signal are the switching period $\Ts$ and the duty cycle $D$ which is the relation between the ``on''-time and the switching period. Given a reference signal $r(t)$ and a carrier signal $s(t)$ the pulsed voltage is generated by{
\begin{equation}
  v_i(t) = \frac{V_i}{2}\,\big(\mathrm{sgn}\left(r(t)-s(t)\right)+1\big), \label{pels:equ:genPWMnat}
\end{equation}
where $\mathrm{sgn}$ denotes the sign function and $V_i$ is the amplitude}. The converter circuit is mathematically described by a system of ordinary or differential-algebraic equations (DAEs), e.g., 
\begin{align}
  \bfA \frac{\dd }{\dd t} \bfx(t) + \bfB\,\bfx(t) &= \bfc(t),\qquad t\in(t_0,T],
  \label{pels:equ:origODEs}
\end{align}

\noindent with given initial value $\bfx(t_0)=\bfx_0$, where $\bfx(t)\in\mathbb{R}^{\Ns}$ is the unknown solution vector consisting for example of currents and voltages, $\bfA$, $\bfB\in\mathbb{R}^{\Ns\times \Ns}$ are matrices, and $\bfc(t)\in\mathbb{R}^{\Ns}$ is the right-hand side containing current and voltage sources, e.g.,  the pulsed voltage $v_i(t)$. The system may be assembled from lumped element descriptions based on loop or (modified) nodal analysis as described in \cite{Estevez-Schwarz_2000aa}. Please note, that we focus on the linear case but the approach can be straight-forwardly generalized, e.g., considering $\bfB=\bfB(\bfx)$. 

\begin{figure}[t]
  \centering
  \includegraphics{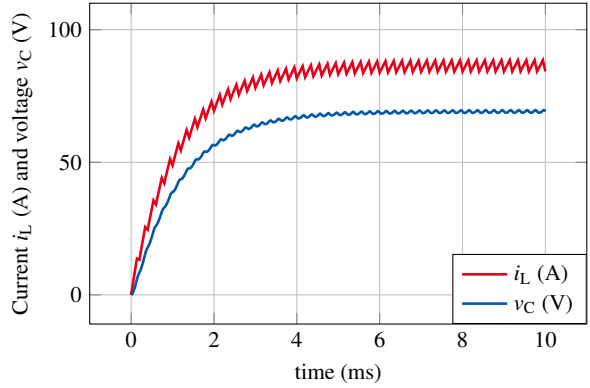}
  \vspace{-0.5em}
  \caption{Exemplary solution of the buck converter depicted in Fig.~\ref{pels:fig:buck_conv}. Switching frequency $\fs=1/\Ts=5\,$kHz.}
  \label{pels:fig:solution}  
\end{figure}
The solution of power converters, e.g., shown in Fig.~\ref{pels:fig:solution}, exhibits the multirate phenomenon: slow variations in the solution require large time intervals until a steady state is reached, i.e., a large end time point $T$ while the fast dynamics due to the switching enforce small time steps. This is the motivation to turn to (parallel) methods that can exploit this multirate behavior. In the following, we focus on the settling process until the steady state is reached. {If one is interested only in the latter, then other methods may also be used, for example the application of Parareal for time-periodic problems is a natural generalization of this work, see, e.g., \cite{Gander_2018aa}}.

\section{Parareal Algorithm}
\label{pels:sec:parareal}

Parareal is an iterative algorithm which is able to accelerate 
the solution of \eqref{pels:equ:origODEs} via parallelization in time. The method originates from \cite{Lions_2001aa} and its 
superlinear convergence is proven in \cite{Gander_2008aa}. The two main ingredients of Parareal are the fine and the coarse propagators. 
We denote by $\mathcal{F}(t,t_0,\bfx_0)$ and $\mathcal{G}(t,t_0,\bfx_0)$ the solutions 
of the {initial value problem (IVP)} \eqref{pels:equ:origODEs} at $t\in(t_0,T]$ obtained with sequential time 
stepping using fine and coarse time steps, respectively.

Partitioning the time interval $t_0=T_0<T_1<\dots<T_N=T$ we write the Parareal iteration: 
for $k=0,1,\dots$ and $n=1,\dots,N$ solve 
  \begin{align}
  \bfX_0^{(k+1)}&=\bfx_0,\label{pels:equ:PR_init}\\
  \bfX_n^{(k+1)}&=\mathcal{F}\big(T_n, T_{n-1},{\bfX^{(k)}_{n-1}}\big)+\mathcal{G}\big(T_n, T_{n-1},{\bfX^{(k+1)}_{n-1}}\big) - \mathcal{G}\big(T_n, T_{n-1},\bfX^{(k)}_{n-1}\big).
  \label{pels:equ:PR}
  \end{align}
The solution operator $\mathcal{F}$ is assumed to deliver a very accurate solution (e.g., using a numerical time-integration method with small time steps $\delta T$) and can be executed in parallel, while $\mathcal{G}$ gives rough information about the solution using a cheap method (e.g., using a numerical method with large time steps $\Delta T_i=T_{i+1}-T_i$) and has to be calculated 
sequentially, cf. \eqref{pels:equ:PR}. 

\medskip

A difficulty in applying Parareal to solve problems with PWM input is that a naive implementation of a coarse propagator using a time-integrator with large time steps will not capture the high-frequency dynamics and may also fail to propagate low-frequency components. A modified Parareal algorithm which still approximately captures the high-frequency behavior was introduced in \cite{Gander_2019aa}.  The idea is to separate the high-frequency (pulsed) components from the low-frequency components, i.e.,
\begin{equation}
  \bfA \frac{\dd }{\dd t} \bfx(t) + \bfB\,\bfx(t) = \underbrace{
  \bar{\bfc}(t)+\tilde{\bfc}(t)}_{=\bfc(t)},
\end{equation}
where $\bar{\bfc}$ can be given as a few low-frequency sinusoids from a (fast) Fourier transform and 
$\tilde{\bfc}(t):=\bfc(t)-\bar{\bfc}(t)$ is the remainder. This allows to define a reduced coarse propagator $\bar{\mathcal{G}}_\textrm{fft}$ which solves
\begin{equation}
  \label{pels:equ:reducedODEs}
  \bfA \frac{\dd }{\dd t} \bfx(t) + \bfB\,\bfx(t) =\bar{\bfc}(t)
\end{equation}
and gives rise to a modified Parareal update formula {with coarse propagator $\bar{\mathcal{G}}_\textrm{fft}$ in \eqref{pels:equ:PR_init}-\eqref{pels:equ:PR}}. 
This modified method converges reliably but possibly with reduced order \cite{Gander_2019aa}. In this paper we propose an alternative method to perform time integration by using the MPDE approach as the coarse propagator.

\section{Multirate PDEs}
\label{pels:sec:mpdes}
\begin{figure}[t]
  \begin{subfigure}{0.5\textwidth}
    \centering
	\includegraphics{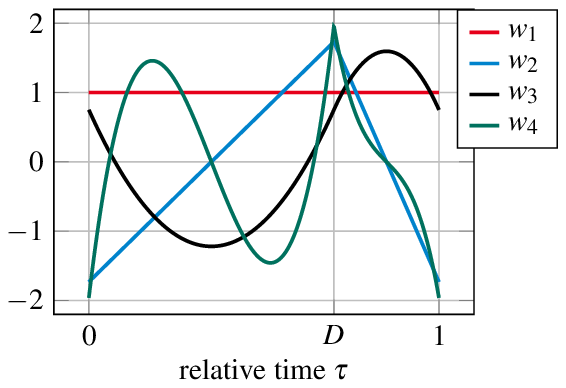}
    \caption{}
  \end{subfigure}
  \begin{subfigure}{0.5\textwidth}
    \centering
	\includegraphics{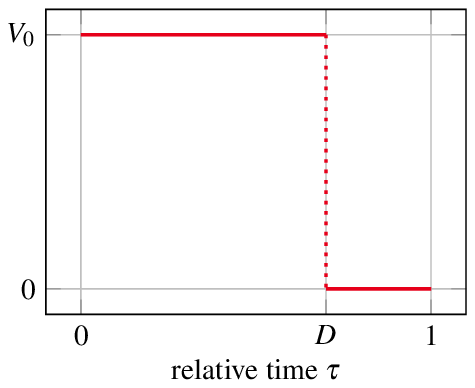}
    \caption{}
  \end{subfigure}
  \caption{
    Construction of basis functions with cusp at relative switching time $D$: (a) PWM basis functions on relative time interval and (b) right-hand side.
    \label{pels:fig:pwmbasis}
  }
\end{figure}

The MPDE approach, which is used for obtaining the coarse solution in Parareal uses the MPDE concept \cite{Brachtendorf_1996aa}.
For the given problem the solution can be conveniently decomposed into a slowly varying envelope and fast periodically varying ripples using the solution expansion \cite{Pels_2019aa}
\begin{equation}
  \xhat_j(t_1, t_2)\doteq\sum_{k=1}^{\Np} y_{j,k}(t_1) w_k(\tau(t_2)) = \bfw^\top\!(\tau(t_2)) \bfy_j(t_1), 
  \label{pels:equ:solExp}
\end{equation}
where $y_{j,k}(t_1)$ are slowly varying coefficients and $w_k(\tau(t_2))$ are a finite set of basis functions ($k=1,\ldots,N_\textrm{p}$) whose periodicity is accounted for by the relative time $\tau(t_2)=\frac{t_2}{\Ts}\text{ mod }1$.
\color{black}
Its application to \eqref{pels:equ:origODEs} yields
\begin{equation}
  \bfA \left(\frac{\partial \bfxhat(t_1,t_2)}{\partial t_1} + \frac{\partial \bfxhat(t_1,t_2)}{\partial t_2}\right) + \bfB\,\bfxhat(t_1, t_2) = \bfchat(t_1,t_2),
  \label{pels:equ:origMPDEs}
\end{equation}
where the relation between the original \eqref{pels:equ:origODEs} and the MPDE \eqref{pels:equ:origMPDEs} solution and right-hand side are given by
\begin{equation}
   \begin{aligned}
     \bfxhat(t, t)=\bfx(t),\quad
     \bfchat(t, t)=\bfc(t).
   \end{aligned}
   \label{pels:equ:relationMPDEsODEs}
\end{equation}
This implies that if a solution to \eqref{pels:equ:origMPDEs} is found, the solution of \eqref{pels:equ:origODEs} can be extracted from it. 
Applying a Galerkin approach along the fast time scale $t_2$ leads to the enlarged equation system
\begin{equation}
  \cA(t_1) \, \frac{\dd  \bfy}{\dd t_1} + \cB(t_1) \, \bfy(t_1) \ = \ \cC(t_1) \, ,
  \label{pels:equ:reducedMPDEs}
\end{equation}
where the matrices are given by \cite{Pels_2019aa}
\begin{align*}
  \cA &=\bfA\otimes\cJ, 
  && \text{with } \qquad
  \cJ = \Ts \int\limits_0^1  \bfw(\tau) \, \bfw^\top \!(\tau)  \,\dd \tau,
  \\
  \cB &=\bfB\otimes\cJ+\bfA\otimes\cQ,
  && \text{with } \qquad
  \cQ = - \int\limits_0^1  \frac{\partial \bfw(\tau)}{\partial \tau } \, \bfw^\top\!(\tau)  \,  \dd \tau,
  \\
  \cC &=\int_{0}^{\Ts}
  \bfchat(t_1,t_2) \otimes \bfw(\tau(t_2))  \, \dd t_2 \, .&&
\end{align*}
Suitable basis functions, which can well represent the ripples in the power converter solution, are, e.g., B-Splines with suitable continuity or the PWM basis functions  \cite{Gyselinck_2013ab}. The latter are global polynomial ansatz functions with $w_1(\tau,D)=1$, $w_2(\tau,D)$ piecewise linear and $w_k(\tau,D)$ is obtained recursively by integrating $w_{k-1}(\tau)$ and orthonormalizing for $3\leq k \leq N_{\textrm{p}}$, see Fig.~\ref{pels:fig:pwmbasis}. It has been shown in \cite{Pels_2019aa} that they are capable of very effectively representing the ripples in linear problems.

\medskip

Finally, equation \eqref{pels:equ:reducedMPDEs} can be time-stepped along $t_1$ by using much larger time steps than are needed to solve \eqref{pels:equ:origODEs} since the fast variations are taken into account by the basis functions. The accuracy of the solution (reconstructed using \eqref{pels:equ:solExp}) increases with $N_{\textrm{p}}$. However increasing $N_{\textrm{p}}$ also makes each time step of an implicit method more costly since an enlarged linear equation system has to be inverted. Nevertheless, even with very few basis functions the reconstructed solution can be expected to capture the main features of the exact solution. This motivates the introduction of another coarse propagator $\bar{\mathcal{G}}_\textrm{mpde}$ in Parareal which solves \eqref{pels:equ:reducedMPDEs} and extracts afterwards the single-time solution according to \eqref{pels:equ:solExp}. 

\section{Numerical Experiments}
\label{pels:sec:numerics}
The proposed approach is applied to the example of the buck converter (see Fig.~\ref{pels:fig:buck_conv}). Its circuit is described by the IVP \eqref{pels:equ:origODEs} given by 
$$
\bfA=
\begin{bmatrix}
L & 0\\
0 & C
\end{bmatrix};
\quad
\bfB=
\begin{bmatrix}
R_{\mathrm{L}} & 1\\
-1 & 1/R 
\end{bmatrix}
\quad
\text{and}
\quad
\bfc(t)=
\begin{bmatrix}
v_i(t)\\
0
\end{bmatrix},
$$
with inductance $L=\num{e-3}\,\SI{}{\henry},$ capacitance $C=\num{e-4}\,\SI{}{\farad},$ resistances $R_{\mathrm{L}}=\num{e-2}\,\SI{}{\ohm}$ and $R=\SI{0.8}{\ohm}.$ The PWM input $v_i(t)$ has the amplitude of $V_i=\SI{100}{\volt}$ and is generated by a {sawtooth carrier signal $s(t) = t \fs\text{ mod }1$} with switching frequency of $\fs=\SI{5}{\kilo\hertz}$ and the reference signal $r(t)=0.7$ {according to \eqref{pels:equ:genPWMnat}}. 
The considered time interval $[0,12]\,$ms is partitioned into $N=40$ windows {for all Parareal variants}. The coarse time step size is {$\Delta T=T/N=\SI{3e-4}{s}$} and the fine propagator uses the time step $\delta T = \num{e-6}\,$s. {All solutions are} obtained with the implicit Euler method. 

First, the classical Parareal method \eqref{pels:equ:PR_init}-\eqref{pels:equ:PR} is applied. {It solves the original system \eqref{pels:equ:origODEs} with the PWM input in both propagators, i.e., $\mathcal{G}$ and $\mathcal{F}$. This} method is compared to two variants where $\mathcal{G}$ is changed to: 1.) $\bar{\mathcal{G}}_\textrm{fft}$ which solves system \eqref{pels:equ:reducedODEs} containing only the DC component instead of the PWM signal on the right-hand side 
(modified Parareal \cite{Gander_2019aa}); 2.) $\bar{\mathcal{G}}_\textrm{mpde}$ which solves \eqref{pels:equ:origODEs} using the MPDE approach with $\Np=1$ and $\Np=3$ with the right-hand-side {$\bfchat(t_1,t_2) = \bfc(t_2)$}.

\begin{figure}[t]
	\begin{subfigure}{0.5\textwidth}
		\centering
		\includegraphics{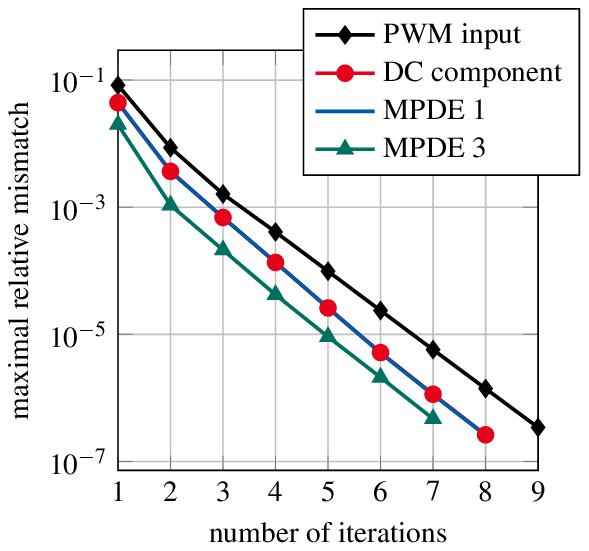}
		\vspace{-1em}
		\caption{\label{pels:fig:convergence}}
	\end{subfigure}
	\begin{subfigure}{0.5\textwidth}
		\centering
		\includegraphics{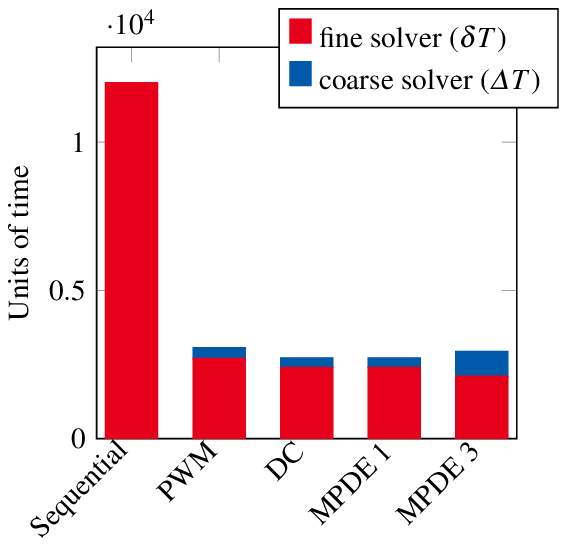}
		\vspace{-1em}
		\caption{\label{pels:fig:cost}}
	\end{subfigure}
	\caption{Convergence of Parareal using different coarse propagators for the buck converter model (a) and units of time for solving the effective linear systems (b).}  
\end{figure}

The maximal relative {mismatch of the solution (`jump') at} the synchronization points $T_n$ for $n=1,\dots,N-1$ is depicted in Fig.~\ref{pels:fig:convergence} for all considered approaches. {The conventional coarse propagator requires always roughly $2$ Parareal iterations more than the MPDE approach with $\Np=3$ to obtain the same accuracy. This is particularly interesting for low accuracy demands, e.g., $10^{-3}$, where we need $4$ vs. $2$ iterations.} The classical Parareal converges up to the relative tolerance of $10^{-6}$ in $9$ iterations. This corresponds to {$2\,700$ and $360$ sequential solutions of linear algebraic systems of size $\Ns=2$ on the fine and the coarse levels, respectively, or $3\,060$ linear systems in total. By the number of sequential solves we mean the number of solver calls which cannot be carried out in parallel (communication costs are neglected).} 
The approaches using the DC component and the MPDE approach with $\Np=1$ both required $8$ iterations ({$2\,400$ fine and $320$ coarse solves, or in total $2\,720$ solutions of linear systems in $2$ variables}). {Finally,} the MPDE approach with $\Np=3$ {basis functions} on the coarse level converged after $7$ iterations, thereby {solving $2\,100$ linear systems of size $\Ns=2$ on the fine level and $280$ linear systems of size $\Ns\times\Np=6$ on the coarse level. }

{For the comparison, let us assume that the overall costs are given by the linear equation system solver which has linear complexity and that solving one linear system of size $\Ns$ requires one unit of time. Then, the classical Parareal takes 
$3060$ units of time. Using MPDE~$1$ (i.e., $\Np=1$) 
or the DC-component as coarse propagator requires only 
$2720$ units of time. Finally, MPDE~3 (i.e., $\Np=3$) takes 
$2940$ units of time. We see that, even in this theoretical setting with an optimal solver, the increased accuracy of the coarse propagator, i.e, application of MPDE~3 with $\Np=3$, does not compensate for the increased effort on the coarse level due to the enlarged equations system, see Fig.~\ref{pels:fig:cost}.
}

Furthermore, from Fig.~\ref{pels:fig:convergence} we see that Parareal with coarse propagator $\bar{\mathcal{G}}_\textrm{mpde}$ using a constant basis function, i.e., $\Np=1$ and the modified Parareal with $\bar{\mathcal{G}}_\textrm{fft}$ using only the DC excitation perform very similarly (if not identically). This resemblance is not surprising since the MPDE~1 approach with $\Np=1$ computes only the envelope of the solution, which is conceptually similar to the modified Parareal with a smooth (in this case constant) coarse input. Finally, we observe that exploitation of more basis functions ($\Np>3$) does not improve the convergence of Parareal, they are similar to the case $\Np=3$.

\section{Conclusions}
\label{pels:conclusions}
In this paper we introduced a novel parallel-in-time algorithm, able to treat systems excited by pulse-width modulated signals. The method extends the two-grid Parareal algorithm by exploiting the MPDE solution approach on the coarse grid. It was applied to the time-domain simulation of a buck converter supplied by a PWM voltage source. Comparison of the proposed algorithm to the standard Parareal method and to the Parareal with reduced coarse dynamics illustrated its faster convergence. Future research will further investigate the similarity of Parareal with the MPDE coarse propagator and the modified Parareal as well as higher order MPDE approaches as coarse propagators.

\begin{acknowledgement}
The authors thank Ruth Vazquez Sabariego from KU Leuven for many fruitful discussions on the MPDE approach.
This research was supported by the Excellence Initiative of the German 
Federal and State Governments and the Graduate School of Computational 
Engineering at Technische Universit\"at Darmstadt, as well as by DFG 
grant SCHO1562/1-2 and BMBF grant 05M2018RDA (PASIROM). 
\end{acknowledgement}

\bibliographystyle{plain}



\end{document}